\documentclass[12pt, a4, reqno,a4paper]{amsart}
\theoremstyle{plain}
\newtheorem*{theorem}{Theorem}

\theoremstyle{remark}

\newcommand{\mg}{\mathrm{SL}_2(\mathbb{Z})}
\newcommand{\Q}{\mathbb Q}

\title{The Ramanujan-Serre differential operators and certain elliptic curves}
\author{Masanobu Kaneko and Yuichi Sakai}
\email{mkaneko@math.kyushu-u.ac.jp, dynamixaxs@gmail.com}

\begin{document}
\maketitle
\begin{abstract} For several congruence subgroups of low levels and their conjugates,
we derive differential equations satisfied by the Eisenstein series of weight 4 and 
relate them to elliptic curves, whose associated new forms of weight 2 constitute
the list of Martin and Ono of new forms given by eta-products/quotients.

\end{abstract}

\section{Introduction}
Let 
\[ E_4(\tau)=1+240\sum_{n=1}^\infty\bigl(\sum_{d\vert n}d^3\bigr)q^n\quad
(q=e^{2\pi i \tau},\, \tau\in\mathfrak{H}: \text{the upper-half plane})
\]
 and 
 \[ E_6(\tau)=1-504\sum_{n=1}^\infty\bigl(\sum_{d\vert n}d^5\bigr)q^n \]
be the standard Eisenstein series of weights 4 and 6 on the modular group 
$\mg$.  The classical relation
\[ E_4(\tau)^3-E_6(\tau)^2=1728 \Delta(\tau), \]
where $\Delta(\tau)=\eta(\tau)^{24}$ with $\eta(\tau)=q^{1/24}\prod_{n=1}^\infty(1-q^n)$ 
the Dedekind eta function, can be seen as the elliptic curve
\begin{equation}\label{cm1}  y^2=x^3-1728 \end{equation}
being parametrized by modular functions 
\[ x=\frac{E_4(\tau)}{\eta(\tau)^8}=\sqrt[3]{j(\tau)},\ y=\frac{E_6(\tau)}{\eta(\tau)^{12}}=\sqrt{j(\tau)-1728}, \]
where $j(\tau)=E_{4}(\tau)^3/\Delta(\tau)$ is the elliptic modular function.

In \cite{guerzhoy2005ramanujan}, Pavel Guerzhoy viewed this parametrization
as a differential equation satisfied by $E_4(\tau)$ and its derivative $\partial(E_4(\tau))=E_6(\tau)$,
where $\partial$ is (a suitable multiple of) the Ramanujan-Serre differential operator, and investigated 
from this viewpoint certain Kummer type congruences satisfied by values
of derivatives of modular forms.
We note that the elliptic curve \eqref{cm1} is isomorphic (over $\Q$) to the minimal curve
\[ y^2=x^3-27 \]
of conductor $36$, and that the associated new form of weight 2 and level $36$
is given by the eta product $\eta(6\tau)^4$, which is seen, as will be justified in the
sequel, to be derived from $\Delta(\tau)=\eta(\tau)^{24}$ as $\eta(6\cdot\tau)^{24/6}$
where the number 6 is the half of the weight of $\Delta(\tau)$.

In the present paper, we show that there exist similar procedures when we consider 
the Ramanujan-Serre differential operator for several congruence subgroups of low levels,
and that the list in the paper \cite{martin1997eta} of Yves Martin and Ken Ono 
of weight $2$ newforms given by eta-products (or quotients) can be recovered.

\section{Main result}

Let $N$ be a positive integer and
\begin{equation*}
\Gamma_0(N)=\biggl\{
				\begin{pmatrix}a&b\\c&d\end{pmatrix}	\in\mathrm{SL}_2(\mathbb{Z})
				\ \biggl|\biggr.\ 
				c\equiv0\bmod{N}
				\biggr\}
\end{equation*}
the  Hecke congruence subgroup of level $N$. 
We start with the eta product $\Delta_N(\tau)$ defined by
\begin{equation*}
\Delta_N(\tau)=\Bigl(\prod_{d | N}\eta(d\tau)\Bigr)^{24/\mu_N},
\end{equation*}
where $\mu_N$ is the index of $\Gamma_0(N)$ in
$\mg$, which is given by 
\[ \mu_N =[\mg\,:\, \Gamma_0(N)]=N\prod_{\substack{p\vert N\\ p: \text{prime}}}\Bigl(1+\frac1p\Bigr). \]
The weight $k_N$ of $\Delta_N(\tau)$ is
\[ k_N= \frac{12\sigma_0(N)}{\mu_N},\quad \sigma_0(N)=\text{ the number of divisors of }N, \]
and the $q$-expansion of $\Delta_N(\tau)$ begins at $q^{h_N}$ where
\[ h_N=\frac{\sigma_1(N)}{\mu_N},\quad \sigma_1(N)=\sum_{d\vert N}d. \]
Using these explicit formulae it is easy to compute all values of $N$ for which the weight $k_N$ is an even integer
and the exponent $h_N$ is an integer. Here is the list (we also give values of $\mu_N$ for convenience): \\

\begin{center}
      \begin{tabular}{c|c|c|c|c|c|c|c|c} 
     $N$   & 1 & 2 & 3 & 5 & 6 & 11 & 14 & 15\\  \hline\hline
    $k_N$     & 12 & 8 & 6 & 4 & 4 & 2 & 2 & 2     \\ \hline
    $h_N$    & 1& 1& 1& 1& 1& 1& 1& 1        \\ \hline
    $\mu_N$ & 1 & 3 & 4 & 6 & 12 & 12 & 24 & 24
            \end{tabular}
   \end{center}
\vspace{1em}
   
\noindent In all cases the exponent $24/\mu_{N}$ in the definition of $\Delta_N(\tau)$ are positive integers.
For these $N$, the form $\Delta_N(\tau)$ is a cusp form of weight $k_N$ on the 
group $\Gamma_0(N)$.  To see this we can use, for instance, the theorem of Honda-Miyawaki \cite[Theorem 1]{HM},
but this is standard anyway and we omit the details here.

Let $P_N(\tau)$ be the logarithmic derivative of $\Delta_N(\tau)$,
\[ P_N(\tau)=q\frac{d}{dq}\log \Delta_N(\tau)=\frac1{\mu_{N}}\sum_{d\vert N}dE_{2}(d\tau), \]
where \[ E_{2}(\tau)=1-24\sum_{n=1}^\infty\bigl(\sum_{d\vert n}d\bigr)q^n \]
is the (``quasimodular'') Eisenstein series of weight 2 on $\mg$.
We introduce the (Ramanujan-Serre) differential operator  $\partial_k^{(N)}$ defined by
\[ \partial_k^{(N)}(f)(\tau)=\frac{k_N}{4}\cdot q\frac{df}{dq}(\tau)-\frac{k}{4}\cdot P_N(\tau)f(\tau)  
\qquad \bigl( q\frac{d}{dq}=\frac1{2\pi \sqrt{-1}}\frac{d}{d\tau}\bigr)\]
for modular forms  $f(\tau)$ of weight $k$. Note that we choose a different normalization
from the standard one by the factor $k_N/4$ (in the classical case of $\mg$ 
most commonly used is  $\partial_k(f)(\tau)=qdf/dq-kE_2f/12$).  From the transformation property of $P_N(\tau)$
(which is quasimodular of weight 2),
we see that if $f$ is modular of weight $k$ on $\Gamma_0(N)$, then
$\partial_k^{(N)}(f)$ is modular of weight $k+2$ on $\Gamma_0(N)$.
(A quick way to see this is to consider the derivative of the weight 0
modular function $f(\tau)^{k_{N}}/\Delta_N(\tau)^{k}$, which is modular of weight 2.)

We further consider the group $\Gamma^\sharp_0(N)$ which is conjugate to $\Gamma_0(N)$
by $\bigl(\begin{smallmatrix} 2 & 1\\ 0 & 2 \end{smallmatrix}\bigr) $:
\begin{equation*}
\Gamma^\sharp_0(N)=\begin{pmatrix}2&1\\0&2\end{pmatrix}^{-1}
					\Gamma_0(N)
				\begin{pmatrix}2&1\\0&2\end{pmatrix}.
\end{equation*}
This group contains $\Gamma_0(4N)$ as a subgroup. Let
\[ \Delta^\sharp_N(\tau)=-\Delta_N(\tau+1/2) \] and 
\[ P_N^\sharp(\tau)=q\frac{d}{dq}\log \Delta_N^\sharp(\tau).\]
Also let $\partial_k^{(N\sharp)}$ be defined by
\[ \partial_k^{(N\sharp)}(f)=\frac{k_N}{4}\cdot q\frac{df}{dq}-\frac{k}4\cdot P_N^\sharp\cdot f  \]
for a form $f$ of weight $k$. 
Here and in the following, we sometimes suppress the variable $\tau$ of modular forms.
Our main result is the following theorem.

\begin{theorem}
i) For $N=1,2,3,5,6$, let $Q_N=Q_N(\tau)$ be any one of the Eisenstein series of weight $4$
on $\Gamma_0(N)$ associated to cusps.  Then the forms $Q_N$ and $\partial_4^{(N)}(Q_N)$
satisfy a homogeneous (with respect to weights) polynomial relation of degree $3$ (and of weight $12$)
over $\Q[\Delta_N]$.
This relation can be written in the form that the pair 
\[ \left(\frac{Q_N}{\Delta_N^{4/k_N}},\frac{\partial_4^{(N)}(Q_N)}{\Delta_N^{6/k_N}}\right) \]
of modular functions satisfies the equation of an elliptic curve $E_N$ over $\Q$. The isomorphism class of $E_{N}$ over $\Q$
is independent of the choice of $Q_{N}$. The conductor of the minimal model of $E_N$ is
$k_N^2N/4$, and the associated new form of weight $2$ is given by
\[ \Delta_N\Bigl(\frac{k_N}2\tau\Bigr)^{2/k_N}. \]

ii) Suppose $N=1,2,5,6$, and let $Q_N$ be as in i). Put $Q_N^\sharp=Q_N(\tau+1/2)$.
Then we obtain a similar differential equation for $Q_N^\sharp$ and an elliptic curve $E_N^\sharp$
whose equation is satisfied by
\[ \left(\frac{Q_N^\sharp}{(\Delta_N^\sharp)^{4/k_N}},
\frac{\partial_4^{(N\sharp)}(Q_N^\sharp)}{(\Delta_N^\sharp)^{6/k_N}}\right). \]
The conductor of the minimal model of $E_N^\sharp$ is $k_N^2 N$ when $N$ is odd and $k_N^2 N/2$ when $N$ is even, and the associated
new form of weight $2$  is given by
\[ \Delta_N^\sharp\Bigl(\frac{k_N}2\tau\Bigr)^{2/k_N}. \]
\end{theorem}

\noindent{\it Remark.}  All powers of $\Delta_N$ and $\Delta_N^\sharp$ appeared in the theorem become
integral powers of products of  etas.\\
 
Before giving the proof, we tabulate in two tables the differential equations, 
the associated elliptic curves (minimal models) and the
new forms. The coefficients of the differential equations depend on normalizations of 
the  Eisenstein series of weight 4. We find that the isomorphism classes of resulting elliptic curves are 
independent, not only of the choices of normalizing constants but also of the choices of the 
Eisenstein series at cusps, and rather surprisingly, we can choose a normalization of each
Eisenstein series to obtain the {\it same} equation.
For $N=1$, there is only one cusp at $i\infty$, and for $N=2,3,5$, there are two inequivalent 
cusps at $i\infty$ and $0$. For $N=6$, there are four cusps,  $1/2$ and $1/3$ being the cusps 
other than $i\infty$ and $0$.  We take the Eisenstein series of weight 4
associated to cusps $i\infty$ and $0$ as
\[
E_{4,N}^{i\infty}(\tau)=c_N
					\sum_{d | N}\frac{\mu(d)}{d^4}
					E_4\Bigl(\frac{N}{d}\tau\Bigr) \]
and 
\[ E_{4,N}^{0}(\tau)=\frac{c_N}{N^{2}}\sum_{d | N}\mu(d)
					E_4(d\tau)\]				
respectively, where $c_N$ is the normalizing constant
\[ c_N=\prod_{\substack{p | N \\ p : \text{prime}}}\Bigl(1-\frac1{p^4}\Bigr)^{-1}\]
and  $\mu(d)$ is the M\"obius function. When $N=6$, as two 
more Eisenstein series $E^{1/2}_{4,6}(\tau)$ and $E^{1/3}_{4,6}(\tau)$ associated to cusps $1/2$ and $1/3$
respectively, we take 
\begin{equation*}
E^{1/2}_{4,6}(\tau)=\frac{3^{2}}{3^4-1}
				\Bigl(E^{i\infty}_{4,2}(\tau)-E^{i\infty}_{4,2}(3\tau)\Bigr)\  \text{ and }\  E^{1/3}_{4,6}(\tau)=\frac{2^{2}}{2^4-1}
				\Bigl(E^{i\infty}_{4,3}(\tau)-E^{i\infty}_{4,3}(2\tau)\Bigr). 
\end{equation*}
\vspace{1em}

The first table below is the list of initial eta products and the differential equations satisfied by Eisenstein series, 
and the second gives minimal models of associated elliptic curves, their conductors, and new forms of weight 2.

\begin{center}\begin{table}[h]
      \begin{tabular}{c|c|c} 
Group & $\Delta_N$ or $\Delta_N^\sharp$ & Differential equation \\ \hline\hline
$\Gamma_0(1)$& $\eta(\tau)^{24}$  & 
$\partial^{(1)}_4(Q_1)^2=Q_1^3-1728\Delta_1$\\[5pt] \hline 
$\Gamma_0(2)$ & $\eta(\tau)^8\eta(2\tau)^8$ & 
$\partial^{(2)}_4(Q_2)^2
=Q_2^3+64\Delta_2Q_2$\\[5pt] \hline
$\Gamma_0(3)$& $\eta(\tau)^6\eta(3\tau)^6$  & 
$\partial^{(3)}_4(Q_3)^2
=Q_3^3+\frac{729}4\Delta_3^2$\\[5pt] \hline
$\Gamma_0(5)$&$\eta(\tau)^4\eta(5\tau)^4$ & 
$\partial^{(5)}_4(Q_5)^2
=Q_5^3-\frac{89}{13}\Delta_5Q_5^2
-\frac{3500}{169}\Delta_5^2Q_5-\frac{125000}{2197}\Delta_5^3$\\[5pt] \hline
$\Gamma_0(6)$& $\eta(\tau)^2\eta(2\tau)^2\eta(3\tau)^2\eta(6\tau)^2$ & 
$\partial^{(6)}_4(Q_6)^2
=Q_6^3-\frac{23}{5}\Delta_6Q_6^2
-\frac{432}{25}\Delta_6^2Q_6
-\frac{1296}{125}\Delta_6^3$\\[5pt] \hline
$\Gamma^\sharp_0(1)$& $\frac{\eta(2\tau)^{72}}{\eta(\tau)^{24}\eta(4\tau)^{24}}$ & 
$\partial^{(1\sharp)}_4(Q_1^\sharp)^2=(Q_1^\sharp)^3+1728\Delta_1^\sharp$\\[5pt] \hline
$\Gamma^\sharp_0(2)$ & $\frac{\eta(2\tau)^{32}}{\eta(\tau)^8\eta(4\tau)^8}$ & 
$\partial^{(2\sharp)}_4(Q_2^\sharp)^2
=(Q_2^\sharp)^3-64\Delta_2^\sharp Q_2^\sharp$\\[5pt] \hline
$\Gamma^\sharp_0(5)$& 
$\frac{\eta(2\tau)^{12}\eta(10\tau)^{12}}{\eta(\tau)^4\eta(4\tau)^4\eta(5\tau)^4\eta(20\tau)^4}$ &
$\partial^{(5\sharp)}_4(Q_5^\sharp)^2
=(Q_5^\sharp)^3+\frac{89}{13}\Delta_5^\sharp(Q_5^\sharp)^2
-\frac{3500}{169}(\Delta_5^\sharp)^2Q_5^\sharp+\frac{125000}{2197}(\Delta_5^\sharp)^3$\\[5pt] \hline
$\Gamma^\sharp_0(6)$& $\frac{\eta(2\tau)^{8}\eta(6\tau)^{8}}{\eta(\tau)^2\eta(3\tau)^2\eta(4\tau)^2\eta(12\tau)^2} $   & 
$\partial^{(6\sharp)}_4(Q_6^\sharp)^2
=(Q_6^\sharp)^3+\frac{23}{5}\Delta_6^\sharp(Q_6^\sharp)^2
-\frac{432}{25}(\Delta_6^\sharp)^2Q_6^\sharp
+\frac{1296}{125}(\Delta_6^\sharp)^3$\\[5pt] \hline 
  \end{tabular}

  \end{table}
   \end{center}

\begin{center}\begin{table}[h]
      \begin{tabular}{c|c|c|c} 
Group & Minimal model of $E_N$ or $E_N^\sharp$ & Conductor & Weight 2 new form\\ \hline\hline
$\Gamma_0(1)$& $y^2=x^3-27$  & 36 &
$ \eta(6\tau)^4$\\[5pt]  \hline
$\Gamma_0(2)$ & $y^2=x^3+4x$ &  32 &
$\eta(4\tau)^2\eta(8\tau)^2$\\[5pt]  \hline
$\Gamma_0(3)$& $y^2+y=x^3$  & 27 &
$\eta(3\tau)^2\eta(9\tau)^2$\\[5pt]  \hline
$\Gamma_0(5)$&$y^2=x^3+x^2-36x-140$ & 20 &
$\eta(2\tau)^2\eta(10\tau)^2$\\[5pt]  \hline
$\Gamma_0(6)$& $y^2=x^3-x-24x-36$ & 24 &
$\eta(2\tau)\eta(4\tau)\eta(6\tau)\eta(12\tau)$\\[5pt]  \hline
$\Gamma^\sharp_0(1)$& $y^2=x^3+27$ & 144 & 
$\frac{\eta(12\tau)^{12}}{\eta(6\tau)^{4}\eta(24\tau)^{4}}$\\[5pt]  \hline
$\Gamma^\sharp_0(2)$ & $y^2=x^3-4x$ & 64 &
$\frac{\eta(8\tau)^{8}}{\eta(4\tau)^2\eta(16\tau)^2}$\\[5pt]  \hline
$\Gamma^\sharp_0(5)$& 
$y^2=x^3-x^2-36x+140$ & 80 &
$\frac{\eta(4\tau)^{6}\eta(20\tau)^{6}}{\eta(2\tau)^2\eta(8\tau)^2\eta(10\tau)^2\eta(40\tau)^2}$\\[5pt]  \hline
$\Gamma^\sharp_0(6)$& $y^2=x^3+x-24x+36$   & 48 &
$\frac{\eta(4\tau)^{4}\eta(12\tau)^{4}}{\eta(2\tau)\eta(6\tau)\eta(8\tau)\eta(24\tau)}$\\[5pt] \hline
  \end{tabular}
 \end{table}  \end{center}

\noindent{\it Proof of Theorem.} 
In each case in i), once a candidate of the differential equation is found, 
it becomes a matter of computing enough Fourier series by computer in order to 
rigorously prove that the equation holds true. This is because each term of the 
differential equations in the table is a holomorphic modular form of weight $12$
on each congruence subgroup, and hence, as is well known (the valence formula
or the Riemann-Roch theorem),  the vanishing of Fourier coefficients up to 
$q^{\mu_N}$ ensures that the function is identically
zero.   Computations of minimal models of elliptic curves
and associated new forms are now standard, and several tables are available
on line.

As for ii), the differential equation for each $\Gamma_0^\sharp(N)$ is only
a translation from that for $\Gamma_0(N)$ by $\tau\mapsto \tau+1/2$.
We should note that the sign of the definition of $\Delta_N^\sharp(\tau)$ is
essential. The computations of datas of elliptic curves are again standard.
This completes the proof of Theorem. \hfill \qed \\

\noindent {\it Remarks.}  1)\  For the remaining cases $N=11,14,15$, the weight of $\Delta_N(\tau)$ is 2 and
this form itself gives the newform associated to an elliptic
curve of conductor $N$. Together with these, the new forms we obtained in Theorem
coincide with the forms in the list of Martin-Ono in \cite{martin1997eta}.  We are very curious to know why this happens, 
what is the reason behind this coincidence.

For $N=11,14,15$, the same procedure as in Theorem 1 does not work for the Eisenstein
series of weight 4. However, we find that certain forms of weight 4 satisfy desired forms of differential
equations.  We briefly describe the form of weight 4 and the differential equation in each case.
The forms $\Delta_{N}$ and $P_{N}$, and the operator $\partial_k^{(N)}$ are as before.
We put for any integer $m$
\[ E_{2,m}(\tau)=\frac{mE_2(m\tau)-E_2(\tau)}{m-1}. \]

For $N=11$, we let 
\[Q_{11}=E_{4,11}^{i\infty}(\tau)-\frac{121}{61} E_{2,11}(\tau) \Delta_{11}(\tau)
-\frac{1274}{915} \Delta_{11}(\tau)^2.\]
Then the equation
\[ \partial^{(11)}_4(Q_{11})^2
=Q_{11}^3-\frac{31}3\Delta_{11}^4 Q_{11}-\frac{2501}{108}\Delta_{11}^6 \]
holds. The minimal model of the elliptic curve $y^2=x^3-\frac{31}3 x-\frac{2501}{108}$
is \[ y^2+y=x^3-x^2-10x-20, \]
which is of conductor 11 and the associated new form of weight 2 is $\Delta_{11}(\tau)=\eta(\tau)^2\eta(11\tau)^2.$

For $N=14$, let 
\[E_{4,14}^{i\infty}(\tau)=\frac{14^4E_4(14\tau)-7^4E_4(7\tau)-2^4E_4(2\tau)+E_4(\tau)}{(7^4-1)(2^4-1)} \]
and set
\[Q_{14}=E_{4,14}^{i\infty}(\tau)-\frac1{900}\bigl(1001E_{2,14}(\tau)-168E_{2,7}(\tau)+73E_{2,2}(\tau)\bigr) \Delta_{14}(\tau). \]
Then we have the equation
\[ \partial^{(14)}_4(Q_{14})^2
=Q_{14}^3-\frac{187}{100}\Delta_{14}^2 Q_{14}^2+\frac{3528}{625}\Delta_{14}^4 Q_{14}-\frac{3863552}{421875}\Delta_{14}^6. \]
The minimal model of the elliptic curve $y^2=x^3-\frac{187}{100}x^2+\frac{3528}{625}x-\frac{3863552}{421875}$
is the conductor 14 curve \[ y^2+xy+y=x^3+4x-6, \]
whose associated new form  is $\Delta_{14}(\tau)=\eta(\tau)\eta(2\tau)\eta(7\tau)\eta(14\tau).$

Finally for $N=15$, let 
\[E_{4,15}^{i\infty}(\tau)=\frac{15^4E_4(15\tau)-5^4E_4(5\tau)-3^4E_4(3\tau)+E_4(\tau)}{(5^4-1)(3^4-1)} \]
and set
\[Q_{15}=E_{4,15}^{i\infty}(\tau)-\frac1{208}\bigl(210E_{2,15}(\tau)-10E_{2,5}(\tau)+9E_{2,3}(\tau)\bigr) \Delta_{15}(\tau). \]
Then we have the equation
\[ \partial^{(15)}_4(Q_{15})^2
=Q_{15}^3-\frac{209}{104}\Delta_{15}^2 Q_{15}^2-\frac{93825}{10816}\Delta_{15}^4 Q_{15}+\frac{860625}{1124864}\Delta_{15}^6. \]
The minimal model of the elliptic curve $y^2=x^3-\frac{209}{104}x^2-\frac{93825}{10816}x+\frac{860625}{1124864}$
is the conductor 15 curve \[ y^2+xy+y=x^3+x^2-10x-10, \]
whose associated new form  is $\Delta_{15}(\tau)=\eta(\tau)\eta(3\tau)\eta(5\tau)\eta(15\tau).$ \\

2)\  In ii) of Theorem, the case of  $\Gamma_0^\sharp(3)$ is missing because  the resulting differential equation is the same as
that for $\Gamma_0(3)$. We have the relation 
\[ \Delta^\sharp_3(3\tau)^{1/3}=
\Delta_3(3\tau)^{1/3}+4\Delta_3(12\tau)^{1/3} \]
between $\Delta^\sharp_3(3\tau)^{1/3}$ and $\Delta_3(3\tau)^{1/3}$.\\

3)\   For $N=1,5,6$, the elliptic curve obtained for $\Gamma_0^\sharp(N)$ is
a quadratic twist of that for $\Gamma_0(N)$, whereas in the case $N=2$, we have a
quartic twist.

\section{Further comments and remarks}

\noindent 1. We have no conceptual reason why we should have started with our eta product $\Delta_N(\tau)$.
In fact, it  is possible that there are  more instances of obtaining elliptic curves
in the way described in our theorem. For example, consider the case of $\Gamma_0(7)$. Let
\[ \Delta_7(\tau)=E_{1,\left(\frac{}{7}\right)}(\tau)\eta(\tau)^3\eta(7\tau)^3\]
be the unique cusp form of weight 4 on $\Gamma_0(7)$
where
\[ E_{1,\left(\frac{}{7}\right)}(\tau)=1+2\sum_{n=1}^\infty\sum_{d|n}\left(\frac{d}7\right)q^n \]
is the Eisenstein series of weight 1 with Legendre character $\left(\frac{}{7}\right)$, and define $P_7$ and $\partial^{(7)}_k$ as
before.  Note that $\Delta_7(\tau)$ has a zero in the upper half-plane and so the 
operator $\partial^{(7)}_k$ does not necessarily send holomorphic modular forms on $\Gamma_0(7)$
to holomorphic ones. But still, we have the differential equation
\[ \partial^{(7)}_4(Q_{7})^2
=Q_{7}^3-\frac{17}{10}\Delta_{7} Q_{7}^2-\frac{637}{100}\Delta_{7}^2Q_7-\frac{45619}{1000}\Delta_7^3 \]
for $Q_7=E_{4,7}^{i\infty}(\tau)$ or $E_{4,7}^0(\tau)$.
The corresponding elliptic curve has the minimal model
\[ y^2=x^3+x^2-7x-52\] of conductor $336=2^4\cdot3\cdot7$.

What is the class of elliptic curves over $\Q$ obtained in the way described in the present paper?
 \\

\noindent 2.  We searched for power series solution
\[ a_0+a_1q+a_2q^2+a_3q^3+\cdots \]
of each differential equation for $\Gamma_0(N)$ in Theorem (the case of $\Gamma_0^\sharp(N)$
being obtained from this by $q\to -q$). When $N=1,2,6$, no solution other than the Eisenstein series 
exists. 

For $N=3$, we find one more solution 
\[ -27\frac{\Delta_3(\tau)}{E_{2,3}^{i\infty}(\tau)}=
-27 q+486 q^2-5103 q^3+43956 q^4-347490 q^5+\cdots, \]
where \[ E_{2,3}^{i\infty}(\tau)=\frac12(3E_2(3\tau)-E_2(\tau)) \]
is the Eisenstein series of weight 2. This solution  is also of weight 4 but has a pole at an elliptic point. 

When $N=5$, we find two more solutions other than the Eisenstein series:
\begin{align*}
&\quad\Bigl(-\frac{18}{13}\pm2i\Bigr)\frac{j_5(\tau)-2\pm11i}{j_5(\tau)+11\mp2i}
\cdot \Delta_5(\tau)\\
&=\Bigl(-\frac{18}{13}\pm2 i\Bigr) q-\Bigl(\frac{32}{13}\pm52 i\Bigr)
   q^2+\Bigl(\frac{2044}{13}\pm384 i\Bigr) q^3-\Bigl(\frac{19696}{13}\pm1256 i\Bigr)
   q^4+\cdots,
   \end{align*}
where $i=\sqrt{-1}$ and 
\[ j_5(\tau)=\eta(\tau)^6/\eta(5\tau)^6=\frac1q-6+9q+10q^2-30q^3+6q^4-25q^5+\cdots \]
is the ``Hauptmodul'' for $\Gamma_0(5)$, which is also given in terms of the forms we have
used as
\[ j_5(\tau)=\frac{E_{4,5}^{i\infty}(\tau)}{\Delta_5(\tau)}-\frac{125}{13}. \]
These solutions are of weight 4
and have poles at elliptic points $\pm2/5+i/5$ where the function $j_5(\tau)$ assumes 
the values 
\[ j_5\bigl(\pm\frac25+\frac{i}5\bigr)=-11\mp2i. \]
The series obtained by dividing by the leading coefficient
\[ \frac{j_5(\tau)-2\pm11i}{j_5(\tau)+11\mp2i}
\cdot \Delta_5(\tau)=q-(17\mp13 i) q^2+(93\mp143i) q^3-(70\mp806 i) q^4+\cdots \]  
have their coefficients in the integer ring $\mathbb{Z}[i]$.  \\
 
3.  It should be able to pursue the congruences studied in the paper \cite{guerzhoy2005ramanujan}
of Guerzhoy (the case $N=1$), at least for $N= 2, 3$ where associated curves also have complex multiplication. We leave
this task for interested readers.

\section*{Acknowledgements}  It is a pleasure to thank Pavel Guerzhoy for explaining the viewpoint 
developed in \cite{guerzhoy2005ramanujan} to the first named author, and for his interest and encouragemant
regarding our work.
\bibliographystyle{amsplain}
\providecommand{\bysame}{\leavevmode\hbox to3em{\hrulefill}\thinspace}
\providecommand{\MR}{\relax\ifhmode\unskip\space\fi MR }
\providecommand{\MRhref}[2]{%
  \href{http://www.ams.org/mathscinet-getitem?mr=#1}{#2}
}
\providecommand{\href}[2]{#2}

\end{document}